\documentclass[article]{elsarticle}

\usepackage{amssymb,amsmath,epsfig,hyperref}

\journal{Journal of Theoretical Biology}
\biboptions{authoryear}

\newcommand*\patchAmsMathEnvironmentForLineno[1]{%
  \expandafter\let\csname old#1\expandafter\endcsname\csname #1\endcsname
  \expandafter\let\csname oldend#1\expandafter\endcsname\csname end#1\endcsname
  \renewenvironment{#1}%
     {\csname old#1\endcsname}%
     {\csname oldend#1\endcsname}}%

\newcommand{\eq}[1]{Eq.~(\ref{eq:#1})}
\newcommand{\fig}[1]{Fig.~\ref{fig:#1}}

\newcommand{\app}[1]{\ref{app:#1}}

\begin{document}

\begin{frontmatter}
\title{\vspace{-20mm}\textbf{Fixation probabilities on superstars, revisited and revised}}

\author[myu]{Alastair Jamieson-Lane\corref{cor1}}
\ead{aja107@math.ubc.ca}

\author[myu]{Christoph Hauert}
\ead{christoph.hauert@math.ubc.ca}

\address[myu]{Mathematics Department, University of British Columbia, 1984 Mathematics Road, Vancouver BC V6T 1Z2, Canada}

\cortext[cor1]{Corresponding author}

\begin{abstract}
Population structures can be crucial determinants of evolutionary processes. For the Moran process on graphs certain structures suppress selective pressure, while others amplify it (Lieberman \emph{et al.} 2005 \emph{Nature} \textbf{433} 312-316). Evolutionary amplifiers suppress random drift and enhance selection. Recently, some results for the most powerful known evolutionary amplifier, the superstar, have been invalidated by a counter example (D\'iaz \emph{et al.} 2013 \emph{Proc. R. Soc. A} \textbf{469} 20130193). Here we correct the original proof and derive improved upper and lower bounds, which indicate that the fixation probability remains close to $1-1/(r^4 H)$ for population size $N\to\infty$ and structural parameter $H\gg1$. This correction resolves the differences between the two aforementioned papers. We also confirm that in the limit $N,H\to\infty$ superstars remain capable of eliminating random drift and hence of providing arbitrarily strong selective advantages to any beneficial mutation. In addition, we investigate the robustness of amplification in superstars and find that it appears to be a fragile phenomenon with respect to changes in the selection or mutation processes.
\end{abstract}

\begin{keyword}
evolution \sep  Moran process \sep evolutionary graph theory \sep structured populations
\end{keyword}

\end{frontmatter}

\section{Introduction}

Populations evolve according to the principles of natural selection and random drift. The balance between the two competing processes is determined by numerous factors, including both population size and structure 
\citep{antal:PRL:2006,burger:Genetics:1994,nowak:Nature:1992b,fu:JSP:2013}. 
The most malignant tumour is unlikely to cause harm if it arises in the outermost layer of skin and is easily brushed aside, and the most imperative model for climate change has limited influence until it has worked its way from a researchers desk, through the literature into policy making and public awareness. Position matters.

One of the simplest and most influential models of stochastic evolutionary processes in finite populations is the Moran process \citep{moran:book:1962,nowak:Nature:2004}. It is based on an unstructured (or well-mixed) population of constant size $N$, where each individual is classified either as a resident (wild type) or mutant. Each type is assigned a constant fitness, which determines its propensity to reproduce. The fitness of wild types is normalized to $1$ and mutants have fitness $r$. An advantageous mutant has $r>1$, a disadvantageous mutant has $r<1$ and a neutral mutant is indistinguishable in terms of fitness, $r=1$. In every time step, an individual is randomly selected for reproduction with a probability proportional to its fitness and produces a clonal offspring that replaces an individual, selected uniformly at random, in the population. This process is repeated until eventually the population has reached one of the homogenous states of all residents, if the mutant went extinct, or all mutants, if the mutant successfully took over the entire population \citep{moran:book:1962,nowak:Nature:2004,lieberman:Nature:2005}. In both cases, the population has reached fixation. In the absence of mutation, the two homogenous states are absorbing.

The Moran process models evolutionary dynamics based on selection and random drift in finite populations: an advantageous mutant has a higher probability, but no guarantee, to reach fixation and, similarly, an inferior mutant is more likely to be eliminated, but not with certainty. The fixation probability of either type is analytically accessible for any given initial configuration. Of particular interest is the fixation probability of a single mutant, $\rho_M$, that arises in an otherwise homogenous population of wild types:
\begin{align}
\label{eq:moran}
\rho_M = &\ \dfrac{1-\frac1r}{1-\frac1{r^N}}.
\end{align}
In the neutral limit, $r\to1$, all individuals in the population are equally likely to end up as the single common ancestor, leading to a fixation probability of $1/N$.

The original Moran process ignores population structures -- but this is easily addressed by arranging individuals of a population on a graph, such that each node refers to one individual and the links to other nodes define its neighbourhood. \citet{maruyama:TPB:1970} and  \citet{slatkin:Evolution:1981} conjectured that the fixation probability of a mutant in this Moran process on graphs remains unaffected by population structures. \citet{lieberman:Nature:2005} proved that this is indeed true for a broad class of structures and, in particular, holds for simple structures such as lattices or regular networks. At the same time, this classification indicated that fixation probabilities, $\rho$, may differ for some structures by tilting the balance between selection and random drift. Evolutionary suppressors enhance random drift and suppress selection ($\rho<\rho_M$ for $r>1$ and $\rho>\rho_M$ for $r<1$), whereas evolutionary amplifiers exhibit the intriguing property to enhance selection and suppress random drift ($\rho>\rho_M$ for $r>1$ and $\rho<\rho_M$ for $r<1$). 

In recent years, various aspects of the Moran process on graphs have been explored, including effects of population structures on fixation probabilities \citep{antal:PRL:2006,broom:PRSA:2008,voorhees:PRSA:2013}, or fixation times \citep{payne:IEEE:2009,frean:PRSB:2013}, as well as computational techniques \citep{shakarian:proceedings:2011,fu:PRE:2009}. However, the most intriguing result remains that even perfect evolutionary amplification appears to be possible: ``The superstar\ldots [has] the amazing property that for large [population sizes] $N$, the fixation probability of an advantageous mutant converges to one, while the fixation probability of a disadvantageous mutant converges to zero.'' \citep{lieberman:Nature:2005}. 

More recently \citet{diaz:PRSA:2013} provided a sophisticated and elaborate counter example that contradicted the fixation probabilities reported in \citet{lieberman:Nature:2005}. Here we identify the problem in the original proof, correct it and report new upper and lower bounds on the fixation probability for superstars. Moreover, for any $r>1$, a graph can be constructed such that $\rho$ is arbitrarily close to $1$, thus confirming the existence of perfect evolutionary amplification.

\section{Model}

\subsection{Moran process on graphs}
\begin{figure}[ht]
\centerline{\includegraphics{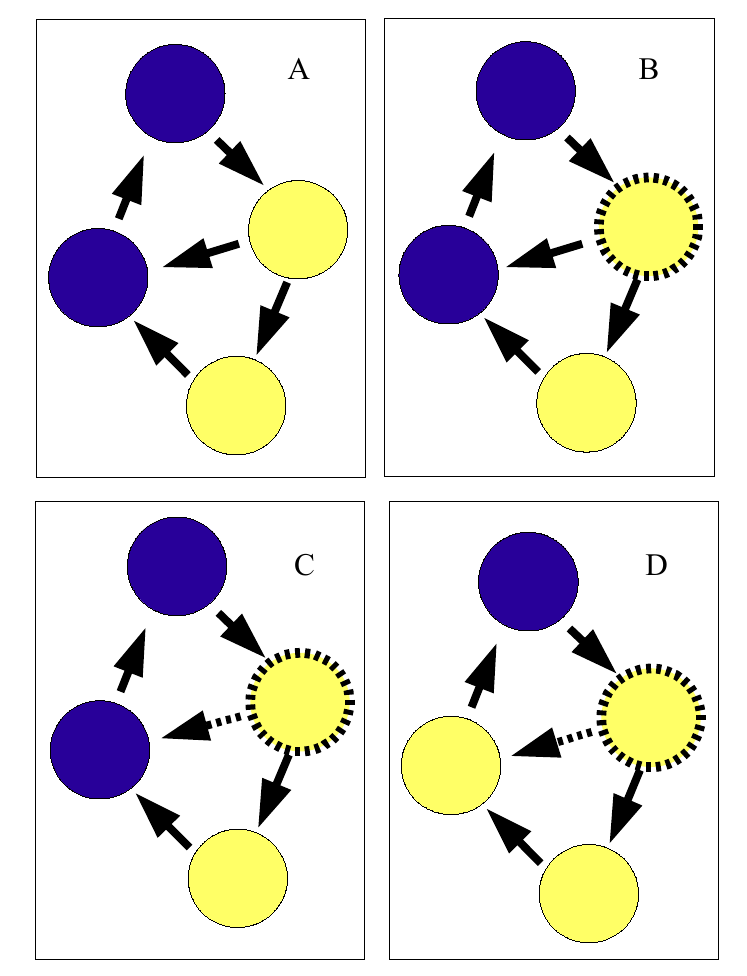}}
\caption{Moran process on a graph. (A) graph structure and distribution of residents (blue) and mutants (yellow). (B) selection: an individual (dashed outline) is selected to reproduce with a probability proportional to its fitness. (C) replacement: a downstream neighbour (dashed arrow) is randomly selected for replacement. (D) reproduction: the neighbour is replaced by the clonal offspring of the upstream reproducing individual.}
\label{fig:invasion}
\end{figure}

Population structure can be represented by assigning individuals to nodes on a graph with links representing each individuals' neighbourhood. The Moran process on graphs follows the same procedure as the original Moran process except for the crucial difference that the offspring does not replace a random member of the entire population but rather replaces a \emph{neighbour} of the reproducing individual, selected uniformly at random (\fig{invasion}).
On directed graphs, the offspring replaces a downstream neighbour by selecting one outgoing link uniformly at random. As before, the population has reached fixation once either one of the absorbing, homogenous states is reached. For any number of mutants, $m$ ($0<m<N$), the fixation probabilities of residents and mutants are both non-zero on strongly connected graphs, i.e. graphs where every node can be reached from any other node through a series of moves between nodes that are connected by links (for directed graphs, only moves in the direction of the link are permitted). If a graph is not strongly connected, then the structure may prevent the spreading or elimination of a mutant type regardless of its fitness and hence the fixation probability for either or both types can be zero.

For the Moran process on graphs, the fixation probabilities are the same as in unstructured populations, c.f. \eq{moran}, provided that the graph is a circulation \citep{lieberman:Nature:2005}. For circulations the sum of weights of all outgoing links is equal to the sum of weights of all incoming links for every node. This means that every node has the same impact on the environment as the environment has on the node.

A graph is an evolutionary suppressor if the fixation probability of an advantageous mutant is less than for the original Moran process, $\rho<\rho_M$. The simplest example is a linear chain: a graph with a single root node, which connects to one end of a (directed) chain of nodes \citep{nowak:PNAS:2003}. Any mutation that does not occur at the root has no chance of reaching fixation. However, if the mutation occurs in the root node it eventually takes over with certainty. Assuming that mutations arise spontaneously and are equally likely in any location, the resulting fixation probability is simply $1/N$, regardless of the mutant's fitness $r$. The linear chain is an example of a graph that is not strongly connected, because the root node cannot be reached from any node in the chain. Evolutionary suppressors are often found when high fidelity copying is of paramount importance, such as in slowing down the somatic evolution of cancer \citep{nowak:PNAS:2003,michor:NRC:2004}.

In contrast, an evolutionary amplifier is a graph, which increases the fixation probability of advantageous mutants as compared to the original Moran process, $\rho>\rho_M$. The simplest evolutionary amplifier is the star graph: a single root node is connected to a reservoir of peripheral leaf nodes through bi-directional links. 
The fixation probability of a single mutant for $N\gg1$ is \citep{lieberman:Nature:2005,broom:PRSA:2008}
\begin{align}
\label{eq:star}
\rho_0 \approx &\ \frac{1-\dfrac1{r^2}}{1-\dfrac1{r^{2N}}}.
\end{align}
On the star, a mutant with fitness $r$ has roughly the same fixation probability as a mutant with fitness $r^2$ in an unstructured population. Thus, the fixation probability of beneficial mutations ($r>1$) is enhanced, but for deleterious mutants ($r<1$) it is reduced. Note that the fixation probability depends on where the single mutant arises. If the mutant is located in the root node then, for $N\gg1$, it is almost certainly replaced in the next time step because one of the $N-1$ reservoir nodes is selected for reproduction. However, if mutants arise at random, then for $N\gg1$ they almost surely arise in the reservoir and the fixation probability is as specified in \eq{star}. Evolutionary amplifiers would seem to provide promising structures for tasks where strong selection is advantageous, such as in the adaptive immune system or in collaboration networks.

\subsection{Superstars}
\begin{figure}[ht]
\centerline{\includegraphics{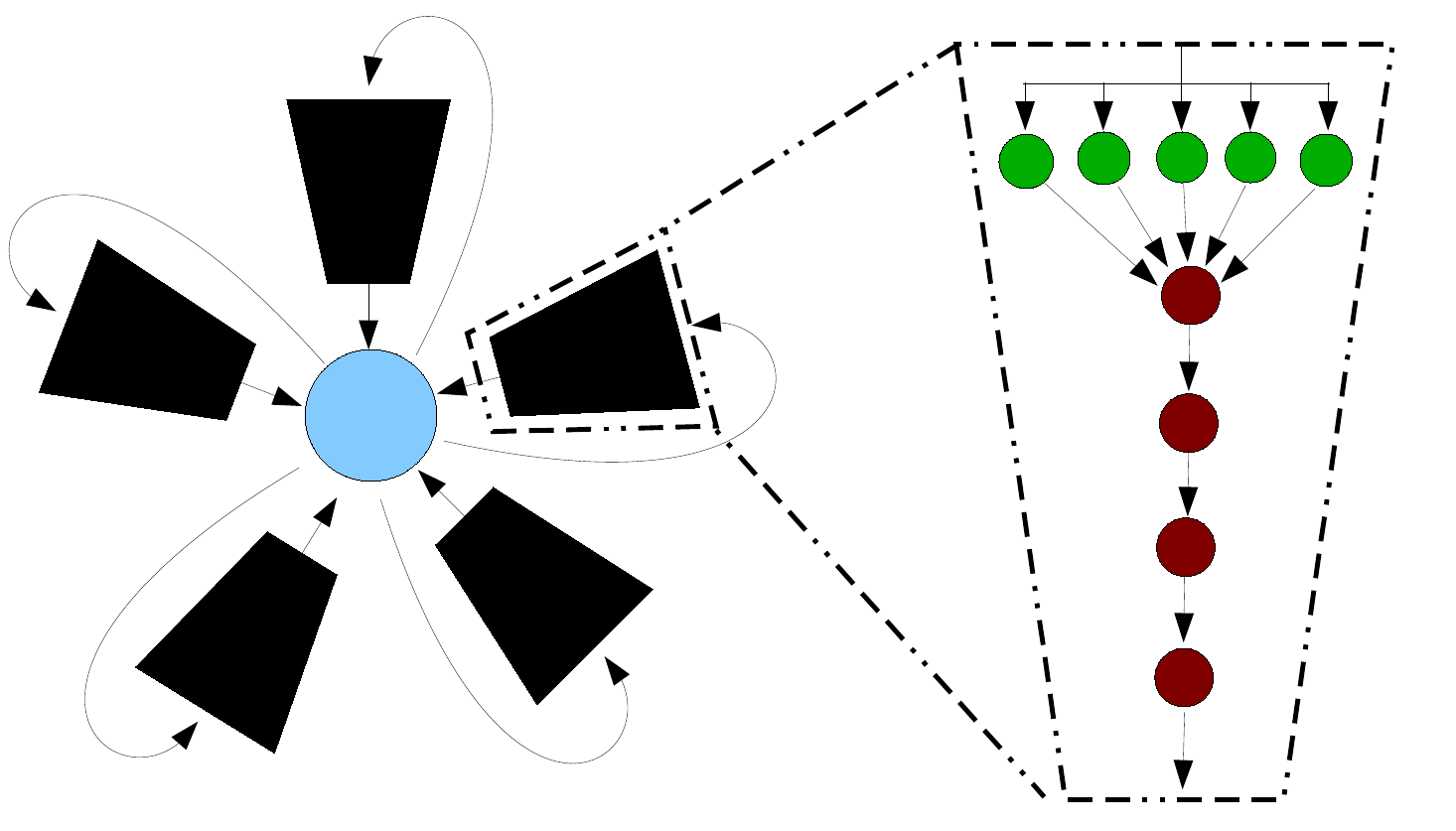}}
 \caption{The superstar consists of three distinct types of nodes: the root node (pale blue), the reservoir nodes (green) and the stem nodes(dark red). The reservoir nodes connect to the start of the stem, the end of the stem connects to the root node and the root node connects to all reservoir nodes in each branch. The depicted superstar has $B=5$ branches each with $L=5$ reservoir nodes and a stem of length $H=4$, which yields a total population size of $N=B(L+H)+1=46$.}
\label{fig:star}
\end{figure}

The two most prominent features of the star graph are the large reservoir where changes occur on a slow time scale, and the bottleneck caused by the hub or root, where changes occur quickly. In particular, the bottleneck introduces a second level for selection to act upon- a mutant needing to reproduce in both the leaves and the hub before it successfully increases its population in the leaves. This basic insight can be exploited to increase evolutionary amplification by elongating the bottleneck and providing further levels where selection can act. Superstars act as a more extreme version of the basic star, and have been proposed as a way to increase evolutionary amplification further \citep{lieberman:Nature:2005}. The superstar consists of a single root node surrounded by $B$ branches (\fig{star}).
Each branch consists of a large reservoir of $L$ nodes feeding into one end of a linear, directed chain of length $H$, the stem. The last stem node in each branch feeds into the root node, which then connects to all reservoir nodes in every branch. The population size is thus given by $N=B(L+H)+1$. Nodes are classified based on their locations on the graph. This classification is designed to simplify discussions but does not affect the rate of reproduction of the individual occupying the node. \citet{lieberman:Nature:2005} report the fixation probability for superstars with $L,B\gg H$ as
\begin{align}
\rho_H  \approx &\ \frac{1-\dfrac1{r^k }}{1-\dfrac1{r^{k N}}},
\label{eq:erez}
\end{align}
where $k=H+2$ is a structural parameter and indicates the number of moves needed to reach any reservoir node from any other reservoir node. This is the number of levels selection can act upon. Consequently it is argued that a single mutant that arises in the reservoir of a superstar with fitness $r$ has approximately the same fixation probability as a mutant with fitness $r^k$ in an unstructured population. This result would then imply that by increasing the length of the stem, the fixation probability, $\rho_H$, of any advantageous mutant, $r>1$, could be brought arbitrarily close to $1$, indicating arbitrarily strong amplification or perfect selection. 

Recently \citet{diaz:PRSA:2013} provided a counter-example demonstrating that the fixation probability in \eq{erez} is too optimistic in the particular case of $H=3$ and thus invalidated the proof in \citet{lieberman:Nature:2005}. In addition, \citet{diaz:PRSA:2013} provide substantial simulation based evidence indicating that \eq{erez} also fails for higher values of $H$.
For the counter-example they show that in the limit $N\to\infty$:
\begin{align}
\rho_3 <&\ 1- \frac{1+r}{2 r^5+r+1}.
\label{eq:diaz}
\end{align} 
This upper bound reflects the probability that a mutant in a reservoir creates a second mutant in any reservoir before getting replaced by resident offspring. Clearly, the fixation probability according to \eq{erez} grows faster with increasing $r$ than \eq{diaz} and for $r>1.42$, results in a contradiction. It turns out that the original proof \citep{lieberman:Nature:2005} was based on an optimistic assumption concerning the amplification along the stem. Taking correlations in the dynamics along the stems into account we obtain new bounds on the fixation probability. More specifically, for $L=B$ we find in the limit $B\to\infty$
\begin{align}
\label{eq:superbounds}
1- \frac{1}{ r^4 (H-1) (1-\frac1r)^2 } \le \rho_H \le 1 - \frac{1}{1+r^4 H}. 
\end{align}
Fixation thus tends to certainty for $H\to\infty$, as suggested by \citet{lieberman:Nature:2005}, while no longer violating the upper bound identified by \citet{diaz:PRSA:2013} for $H=3$. 

In the following we borrow a number of valuable concepts and techniques from both articles, adding and extending where necessary. Exact bounds on the error terms for finite populations are provided in the appendices.

\section{Proof}
\label{sect:proof}%
For the proof of the fixation probability on superstars, we follow the tradition of \citet{lieberman:Nature:2005,diaz:PRSA:2013} and consider superstars with many branches, $B$, and large reservoirs, $L$. More specifically, we study the dynamics of a single branch in detail, and use this to determine the much slower dynamics of changes in the reservoirs.
For any given stem length, $H>2$, the following arguments become exact in the limit $B,L\to\infty$. In practice, we obtain good approximations for $H\ll B, L$. In an effort to reduce notation and increase clarity we assume in the following that $H \ll B=L \approx \sqrt{N}$. However, in the full proof (see appendix) we only require that $B$ and $L$ scale with $N$ but not that they are of the same size.

If mutations arise spontaneously and with equal probability at any node, then the initial mutant almost certainly arises in a reservoir node, because reservoir nodes vastly outnumber nodes of all other types. This marks the starting point for the remainder of our proof (for details see \app{timeZero}).

On occasion, we need to refer to the total fitness of our superstar population at a given time, $F_t$, with $N<F_t<rN$. However, all instances of $F_t$ cancel throughout the proof and hence we do not need to keep track of its exact value. Moreover, various necessary approximations introduce different error terms that are accounted for in full detail in the appendices. All error terms tend to zero as $B,L \to \infty$. It is sufficient to assume that $B=L$ as we take this limit, however other relations are also possible. The exact restrictions on how we can take these limits can be found in \app{summary}, but for now let us simply assume that limits are taken simultaneous, with some suitable relation binding $B$ and $L$ togeather.

\subsection{\label{sect:timescale}Timescales}
Different nodes get updated at different rates. More precisely, any given node is updated if one of its upstream neighbours reproduces and the node of interest is chosen for replacement. Here we follow the convention that all nodes have weighted out-degree of $1$, and that weight is evenly distributed among outgoing links (each link has weight $1/d$ for a node with $d$ outgoing links). Thus, nodes with high (weighted) in-degrees, tend to be short lived, while nodes with low (weighted) in-degrees, tend to be long lived.

Assuming $r\ll N$, every node is selected for reproduction approximately with probability $1/N$. The root node has an in-degree of $B$ and all its upstream neighbours have out-degrees of $1$, hence it updates with a probability close to $B/N\approx 1/\sqrt{N}$. Recall that we assume $H \ll B=L \approx \sqrt{N}$. Similarly, reservoir nodes are replaced with probability of approximately $1/N^2$, the first stem node with probability on the order of $1/\sqrt{N}$, and all other stem nodes with probability of approximately $1/N$. For $N\gg1$, this results in three different timescales: the slowest for reservoir nodes that get replaced, on average, only once in $N^2$ time steps; an intermediate timescale for the stem nodes (with the exception of the  first node), which get replaced once in $N$ time steps; and a fast timescale for the root node as well as the first stem node in each branch, which update once in $\sqrt{N}$ time steps, respectively.

For $N\gg1$, it is possible to separate the three timescales and analyze the dynamics of the different types of nodes individually. More specifically, this allows us to focus on the intermediate timescale associated with the dynamics in the stem, while treating the state of the slowly updating reservoir nodes as constant and the fast updating nodes as random variables. 
In the following, we derive the evolutionary dynamics for the top, middle and bottom of the stem in a single branch. The results determine the slow dynamics of reservoir nodes and describes the early stages of the invasion process, when mutants are rare among the reservoir nodes. This allows us to derive upper and lower bounds on the fixation probabilities.

\subsection{Top of stem}
The first node of the stem gets replaced on the fast time scale, which allows us to treat its state as a random variable. Initially, out of all $L$ upstream neighbours in the reservoir of the corresponding branch, only one is a mutant. Hence, at any given time step, the top node is occupied by a mutant with probability close to $r/L$. This mutant reproduces with a probability $r/F_t$ and hence the probability that a mutant is placed in the second stem node is approximately $r^2/(F_t L)$ in each time step (for error terms, see \app{timeZero}).

\subsection{Middle of stem}
The structure of the stem causes the state of any given stem node to be highly correlated with its neighbours, both upstream and downstream. More specifically, if a mutant reproduces it is highly likely to end up replacing its own offspring. 
This correlation had been neglected in \citet{lieberman:Nature:2005}. In the absence of correlations, whenever a mutant in the stem reproduces, it almost certainly replaces a resident since residents are more common. Along a stem of length $H$ this results in an overall amplification of $r^H$. However, due to correlations, if a mutant in the stem reproduces, it \emph{likely} replaces its former offspring and hence diminishes the resulting amplification.

\begin{figure}[ht]
\centerline{\includegraphics{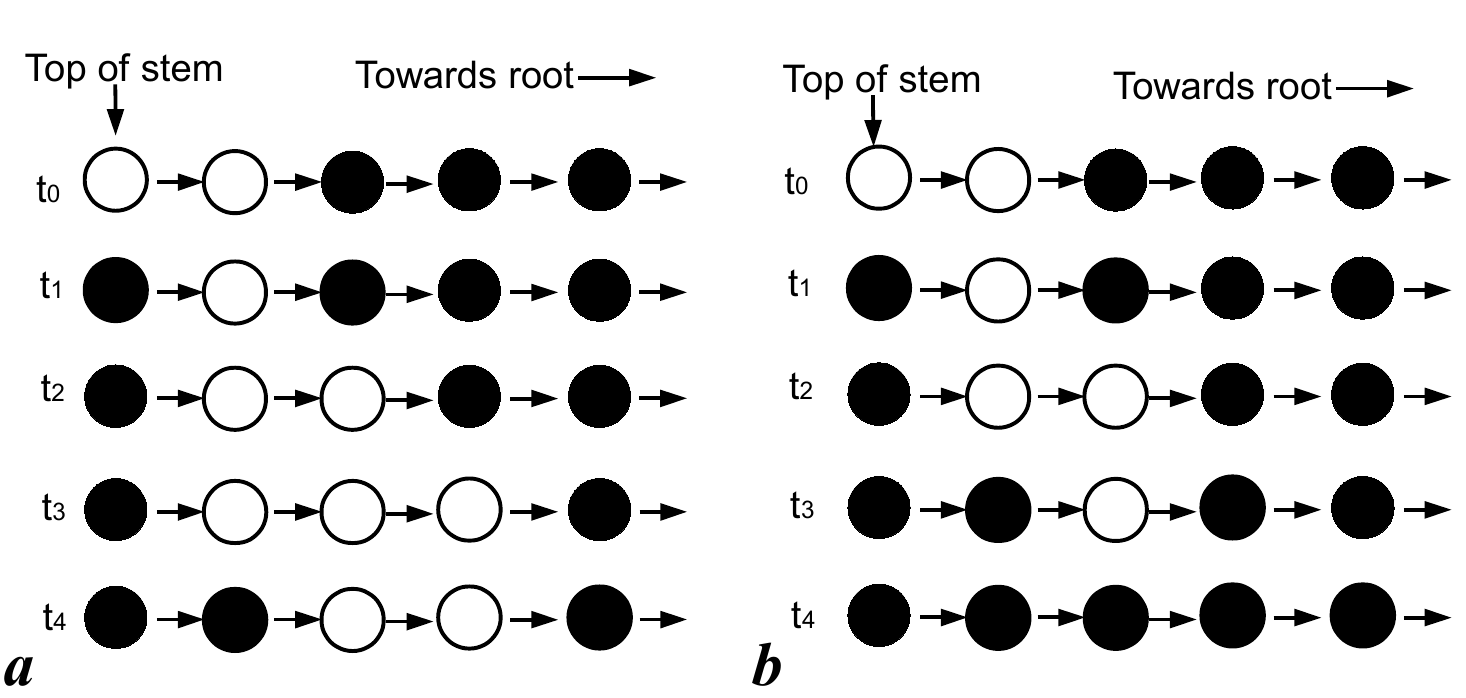}}
 \caption{Two possible histories of a train of mutants (white) proceeding along a stem filled with residents (black).
\textbf{\textsf{a}} We begin with two mutants ($t_0$). The top node is quickly replaced by a resident (on the fast time scale) ($t_1$). Some time later the remaining mutant reproduces ($t_2$), and then the new top node reproduces again ($t_3$). Finally we lose a single mutant from the back of out train ($t_4$). This general growing pattern applies whenever $r>1$. 
\textbf{\textsf{b}} We begin with two mutants ($t_0$), and immediately lose the back mutant of our train ($t_1$). The front of the train reproduces, creating a second mutant ($t_2$), but both fall prey to bad fortune (or low fitness) and are removed ($t_3, t_4$). This behaviour is likely when $r<1$, but even for beneficial mutations many trains do not reach the end of the stem. 
}
\label{fig:trains}
\end{figure}

Simulations nicely illustrate the characteristic features of the dynamics along the body of the stem: clusters of mutants begin at the top of the stem, then grow and move along the stem. In the following, we refer to these clusters as \emph{trains}. A train moves forward and increases in length whenever the front mutant reproduces, which happens at a rate $r/F_t$, but shrinks whenever a resident reproduces and replaces the back end of the train, which occurs at a rate $1/F_t$, see \fig{trains}. 
Thus, as the train moves along the stem, the train length for beneficial mutants increases, on average.

Note that for small superstars with a single node in the stem body, which corresponds to $H=2$ (or $k=4$), the two stem nodes are indeed uncorrelated. However, for $H>2$ this assumption breaks down and results in an overestimation of the fixation probabilities as pointed out by \citet{diaz:PRSA:2013}.

In order to link the stem dynamics to the slow timescale of reservoir nodes, we need to know the expected train length, $T$, by the time the train first reaches the root end of the stem.
The history of any given train can be represented as a random sequence of increments and decrements with a bias that increments are $r$ times more likely.
Essentially, we need to sum over all possible sequences of increments and decrements given an initial train length of $1$ and weigh the resulting train length with the probability of the respective sequence. All sequences where trains disappear (zero train length) are omitted (zero weight). We can count these sequences that result in extinction using the reflection principle (see \app{Mirror}). This method yields the expected train length:
\begin{align}
T =&\ \left(\frac{r}{1+r}\right)^{\!\!H-2}\ \sum_{z=1}^{H-1} (H-z)\!\left(\frac1{1+r}\right)^{z-1} \left[{H-4+z \choose z-1} - {H-4+z \choose z-2}\right].
\label{eq:trains}
\end{align}
For $H\geq2, r>1$ simple bounds for $T$ exist (see \app{simpleT}):
\begin{align}
(H-1) \left(1- \frac1r\right)^2 \le T \le H.
\label{eq:trainsasym}
\end{align}
This indicates that for $r>1$ the train length $T$ grows approximately linearly with increasing stem length $H$.

\subsection{Bottom of stem}
Whenever a train reaches the root end of the stem, its mutants compete with the resident nodes from the other branches to occupy the root node. Since the root node is updated on the fast timescale we can again treat its state as a random variable. Thus, once a train has reached the root end of the stem, the root node is a mutant with probability close to $r/B$. Thus, as long as the train sits at the root end of the stem, the probability in any given time step that the root node is a mutant, reproduces and creates a second mutant in any reservoir is $r^2/(F_t B)$. However, the train is simultaneously eroded from behind, with train mutants being replaced by residents with probability $1/F_t$. Thus, the train remains at the root end for $T F_t$ time steps, on average. Put together, this means that any given train succeeds in producing a second mutant in any reservoir with a probability close to $ r^2 T/ B$ (for detailed error bounds see \app{stemBase}).

\subsection{Slow dynamics in reservoirs}
At any given time step, the probability of losing the initial mutant in the reservoir is $1/(F_t BL)$. Based on the dynamics in the stem, we derive the per time step probability that a second mutant is generated in the reservoir of any branch as the product of the probability that a train is generated and the probability that the train succeeds in producing a second mutant, which yields approximately $r^4 T/(F_t B L)$.  
Thus, the probability to eventually go from one to two mutants in the reservoirs, as opposed to losing the initial mutant, is close to
\begin{align}
\frac{r^4 T}{1+r^4 T}.
\label{eq:1to2}
\end{align}
Since $T$ can be made arbitrarily large (by increasing the stem length $H$, see \eq{trainsasym}), the transition from one to two mutants becomes almost certain and, similarly, the probability of losing the initial mutant becomes vanishingly small.

\subsection{Upper bound on fixation probability}
To find an upper bound on the fixation probability, $\rho_H$, we note that before our mutant can reach fixation, the superstar must first transition from a state with one mutant in a reservoir to a state with two mutants in the reservoirs. Thus, an upper bound on this transition probability serves as an upper bound on the mutant fixation probability. Moreover, the upper bound can be made independent of $T$ by assuming that all trains have the maximum possible train length.
Thus, in the limit of large $B$ and $L$ we find
\begin{align}
\rho_H \le 1-  \frac{1}{1+r^4 T} \le 1-  \frac{1}{1+r^4 H}.
\label{eq:upperBound}
\end{align}
For any given $H, r$ we can find $T$ explicitly using \eq{trains}. In particular, we note that for $H=3$, we find $T= 2r/(r+1)$, thus recovering the upper bound identified in \citet{diaz:PRSA:2013}.

\subsection{Lower bound on fixation probability}
We find a lower bound on the fixation probability by approximating the dynamics of the system with a random walk. This random walk has a forward bias given by \eq{1to2} as long as mutants are rare, and we assume no forward bias otherwise. Because even for larger numbers of mutants the forward bias persists (but there is no simple way to quantify the bias) we obtain a lower bound of the fixation probability, $\rho_H$.

For any finite number of steps, a sufficiently strong initial bias would suffice to ensure that the random walk eventually reaches fixation with high probability. However, the limit $N\to\infty$ also requires an arbitrarily large number of forward steps. In order to resolve the interplay between these two limiting behaviours we set up a martingale and apply the optional stopping theorem \citep[p.~210]{klenke:book:2006} (see \app{LowerFix} for details). In the limit of large $B$ and $L$ we find:
\begin{align}
\rho_H \ge \hat \rho_H = 1- \frac{1}{ r^4 T}   \ge 1- \frac{1}{ r^4 (H-1) (1-\frac1r)^2}.
\label{eq:rhoH}
\end{align}
Once again we note that for any given $H, r$ we can find $T$ explicitly, \eq{trains}. Combined with the upper bound, \eq{upperBound}, this means that $\rho_H$ must exist in a narrow window, \eq{superbounds}.

\section{Robustness}
Unfortunately, it turns out that the intriguing feature of superstars to act as evolutionary amplifiers holds only under very specific conditions. Here we discuss the most important requirements. 

\subsection{Selection \& Sequence}
The original Moran process is formulated as a fecundity based birth-death process, that is, fitness affects the rate of birth (reproduction) whereas death (replacement) occurs uniformly at random. Alternatively, fitness could just as well affect survival such that birth events occur uniformly at random but death events might, for example, occur with probability inversely proportional to fitness. Similarly, the sequence of events could be reversed such that first an individual dies and then the remaining individuals compete to repopulate the vacant site. This yields a total of four distinct scenarios: $Bd$, $bD$, $dB$ and $Db$, where capital letters refer to the fitness dependent selection step. The original Moran process corresponds to $Bd$ and the fixation probability is given in \eq{moran}. In unstructured populations the four dynamical scenarios result in only marginal differences in fixation probabilities. However, they can have crucial effects on the evolutionary outcome in structured populations \citep{ohtsuki:Nature:2006,ohtsuki:JTB:2006b,zukewich:PlosOne:2013}. \citet{frean:preprint:2008} examine all four cases for both complete graphs and star graphs, showing that stars act as evolutionary suppressors in both the $dB$ and $Db$ cases, and are significantly less effective in the $bD$ case compared to the original $Bd$ case. Similar results apply to superstars:

\paragraph{$\mathbf{bD}$ updates:} 
For the birth-death process with selection on survival, mutants only gain any advantage whenever the root node reproduces. Whenever any other node reproduces, there is only a single downstream node, and thus no opportunity for competition, rendering any fitness advantage irrelevant.
This lack of advantage in the stem leads to an expected train length of $1$, regardless of stem length or mutant's fitness. The chance of launching a successful train in a given time step is $1/(NLB)$ and the chance of replacing the original mutant is $1/(NBLr-N(r-1))$. This results in a bias of approximately $\frac{r}{1+r}$ for the initial mutant to eventually create a second reservoir mutant -- the same bias as for the original Moran process. Thus, we might expect fixation probabilities similar to the original Moran process on $BL$ nodes, and certainly nowhere close to the amplification observed for $Bd$ updates.

\paragraph{$\mathbf{Db}$ updates:} 
For the death-birth process with selection on survival, the prospects of mutants drop even further. The probability to successfully place even a single offspring in the top of the stem is only $r/(L+r)$. Note that for death-birth processes the top of the stem no longer changes on the fast timescale and hence trains start at the top instead of the second node. As the train propagates along the stem it tends to grow because the mutant at the back of the train is less likely to die than the resident in front of it, leading to the same train dynamics observed for the $Bd$ process. 
Upon reaching the end of the stem the train competes with the other branches for control over the root node and succeeds with probability near $rT/B$ (over the lifetime of the train). Once a mutant occupies the root, it is predestined to have many offspring -- in each time step a reservoir node dies with high probability and gets replaced by an offspring of the mutant in the root node, whereas the probability is low that the root node is replaced. More specifically, we expect $rN/(1+r)$ reservoir nodes to become mutants before the root node is replaced. At that point it is reasonable to assume that mutants reach fixation with high probability. We conjecture that the probability of mutant fixation on the superstar is close to the probability of a mutant eventually being placed in the root node. Thus, we expect a fixation probability close to $r^2T/BL$. This result is significantly \emph{less} than the almost certain fixation in the original Moran process (c.f. \eq{moran}). The result does, however, match well with the $1/N$ scaling found for the fixation on stars \citep{frean:preprint:2008}.

\paragraph{$\mathbf{dB}$ updates:} 
The final case is the death-birth process with selection on reproduction. Once again the probability of placing a mutant offspring in the stem before losing the reservoir mutant is near $r/L$, but now without further benefits along the stem. Consequently, trains that do reach the root still have an expected length of $1$. Thus, each train has a probability of roughly $r/B$ for claiming the root node, which then produces $N/2$ mutants in the reservoirs, on average -- enough to suggest fixation with high probability, but less than for $Db$. Thus, we conjecture fixation probabilities near $r^2/N$ -- the worst outcome of the four scenarios. Once again, we note the significant penalty as compared to the original Moran process as well as the similarities to the $1/N$ scaling for stars \citep{frean:preprint:2008}.

\subsection{Mutations}
Even though we did not explicitly model the process of mutation, we implicitly assumed that mutations are rare and arise spontaneously in any node selected uniformly at random. For the superstar this means that most mutations arise in a reservoir node -- simply because the overwhelming majority of nodes are reservoir nodes.

An alternative and equally natural assumption is that mutations arise during reproduction events. Such a change does not affect the fixation probabilities in the original Moran process. However, in highly heterogenous population structures crucial differences in the fixation probabilities can arise because mutants preferentially arise in certain locations \citep{maciejewski:PLoSCB:2014}. For superstars, when using the $Bd$ or $bD$ update rules, mutants most likely arise at the top of a stem. This is an unfortunate position because the mutant is highly likely replaced before reproducing even once -- extinction is almost certain. In contrast, for $Db$ and $dB$, mutants again most likely arise among the reservoir nodes -- but for those updates superstars do not act as evolutionary amplifiers.

Even though the dynamical properties of superstars are intriguing, the list of caveats demonstrates that the evolutionary amplification is highly sensitive to the details of the model -- maybe this is the reason that superstar-like structures have not been reported in nature.

\section{Conclusion}
Superstars represent the most prominent representatives of evolutionary amplifiers -- structures that are capable of increasing selection and suppressing random drift. For $r>1$ and in the limit of large $N$ we have derived upper and lower bounds for the fixation probability, $\rho_H$:
\begin{align}
\label{eq:superfinal}
1- \frac{1}{ r^4 T} \le \rho_H \le 1 - \frac{1}{1+r^4 T},
\end{align}
where $(H-1)(1-1/r)^2 \le T \le H$.

Even though fixation probabilities can be made arbitrarily close to $1$ on large superstars and sufficiently large $H$, the fixation probability remains bounded away from $1$ for any finite graph. As a concrete example, consider $r=2$ and $H=50$, which yeilds $T\approx 13.25$ and $0.995283 \le \rho_{50} \le  0.995306$ in the limit of large $N$. Similarly, a sizeable, finite superstar with $B=L=5000$ ($N\approx2.5\cdot10^5$) yields $0.985323\le \rho^N_{50} \le0.995375$, which includes all error terms (see \app{summary}). In contrast, the fixation probability for a similarly sized isothermal graph (e.g. a lattice, complete or random regular graph) is just short of $0.5$.

The upper bound for $\rho_H$ in \eq{superfinal} results in a contradiction with the originally reported fixation probability, \eq{erez}, for sufficiently large $r$. For the specific case of $H=3$ the discrepancy was pointed out in \citet{diaz:PRSA:2013}. At the same time, the lower bound for $\rho_H$ in \eq{superfinal} confirms that superstars are indeed capable of providing an arbitrarily strong evolutionary advantage to any beneficial mutation, as suggested in \citet{lieberman:Nature:2005}. Using symmetry arguments, it also follows that for $r<1$ the fixation probability can be made arbitrarily small, as required for a perfect evolutionary amplifier (see \app{briefNote}).

In the case $H=2$ (or $k=4$ in \citet{lieberman:Nature:2005}) we obtain an expected train length of $T=1$ and recover the original bias, $r^4/(1+r^4)$. Discrepancies arise only for $H\geq3$ (or $k\geq5$) but those cases were not included in the simulations in \citet{lieberman:Nature:2005}. For $H=3$, we obtain $T=2r/(1+r)$, which results in a bias of $2r^5/(1+r+2r^5)$ and recovers the upper bound reported by \citet{diaz:PRSA:2013}. Extending the technique in \citet{diaz:PRSA:2013} to higher values of $H$ numerically, we find that the upper bounds found match \eq{superfinal}. 

An appropriately skeptical reader might ask why the theory presented here should be trusted over those previously presented in the literature -- after all, both claim to offer rigorous proof. First, we note the agreement between predictions made here, and both \citet{lieberman:Nature:2005} and \citet{diaz:PRSA:2013} for the appropriate values of $H$. Second, we identify correlations between neighbouring stem nodes as the cause for the discrepancies between the two previous papers. Finally, we invite readers to scrutinize the proof offered here most thoroughly. Superstars have already presented unexpected subtleties, and as always, we need caution and vigilance to discern between scientific selection and random drift.

\subsubsection*{Acknowledgments}
This work was supported by the Natural Sciences and Engineering Research Council of Canada (NSERC) and the Foundational Questions in Evolutionary Biology Fund (FQEB), grant RFP-12-10.


\bibliographystyle{prsb}
\bibliography{ET}

\appendix

\section{\label{app:timeZero}Initial conditions}

If mutations arise spontaneously and with equal probability in any node then the initial mutant arises in a reservoir node with probability $BL/(BL+1+HB)$. This probability can be made arbitrarily close to one, for suitably large $L$ or $B$. The mutant arises in a stem or root node with probability 
\begin{align}
\label{eq:ep0}
\epsilon_0 = \frac{1+HB}{BL+1+HB}.
\end{align} 
Thus, the final fixation probability, $\rho_H$ (see \eq{superfinal}), needs to be multiplied by $1-\epsilon_0$ to account for this possibility. Because $\epsilon_0$ is only used to derive the lower bound on $\rho_H$, assuming extinction of all mutants not arising in a reservoir node preserves inequalities by effectively ignoring the small possibility that an invading mutant placed in stem nodes or the root node could still reach fixation.

The first node of the stem gets replaced on the fast time scale, which allows us to treat its state as a random variable. However, at early stages of invasion, only one of the $L$ upstream neighbours is a mutant. Hence, at any given time step, the top node is occupied by a mutant with probability 
\begin{align}
\frac{r}{L-1+r}= \frac{r}{L} (1-\epsilon_1)\\
\label{eq:ep1}
\epsilon_1= \frac{r-1}{L+r-1}
\end{align} 
This mutant reproduces with a probability $r/F_t$ and hence the probability that a mutant is placed in the second stem node is $r^2/(F_t (L+r-1))$ in each time step. However, we need to account for the possibility that the initial mutant in the reservoir is replaced \emph{before} the first node in that stem. On a given time step, the chance that the reservoir mutant is replaced by a resident is less than $1/(F_t BL)$. Conversely, the probability of the first node in the chain being replaced exceeds $L/F_t$. Thus the chance that the initial mutant gets replaced before its offspring replaces the first node in the chain is 
\begin{align}
\label{eq:ep2}
\epsilon_2 < 1/(1+ BL^2).
\end{align} 
For our proof we assume that the first node in any chain can be treated as a random variable. The above error term accounts for any slight discrepancy caused by our initial conditions.

\section{\label{app:Mirror}Expected train length $T$.}
Mutants placed in the main body of the stem (excluding the first node, which updates on a fast timescale) propagate down the stem in trains. Trains grow at one end as mutants reproduce, and shrink at the other end as mutants are replaced by residents (see \fig{trains}). At any time, $t$, the state of a train is given by two integers: $A_t$ and $Z_t$. Here $A_t$ refers to the position of the mutant at the front of the train, and $Z_t$ refers to the position of the resident directly behind the train. The current length of the train is thus given by $A_t-Z_t$. Because in most time steps no change occurs in this particular stem, we consider a condensed process, which only accounts for events that change the state of the train. This means that $A_t$ increases with probability $r/(1+r)$ while $Z_t$ increases with probability $1/(1+r)$. Thus, for beneficial mutants the length tends to increase as the train progresses down the stem. If at any time $Z_t\geq A_t$ the train has vanished and the stem is cleared of mutants. In this case, we say that the train, which ``arrives'' at the end of the stem, has length zero.

\begin{figure}[th]
\centering
\includegraphics[width=\textwidth]{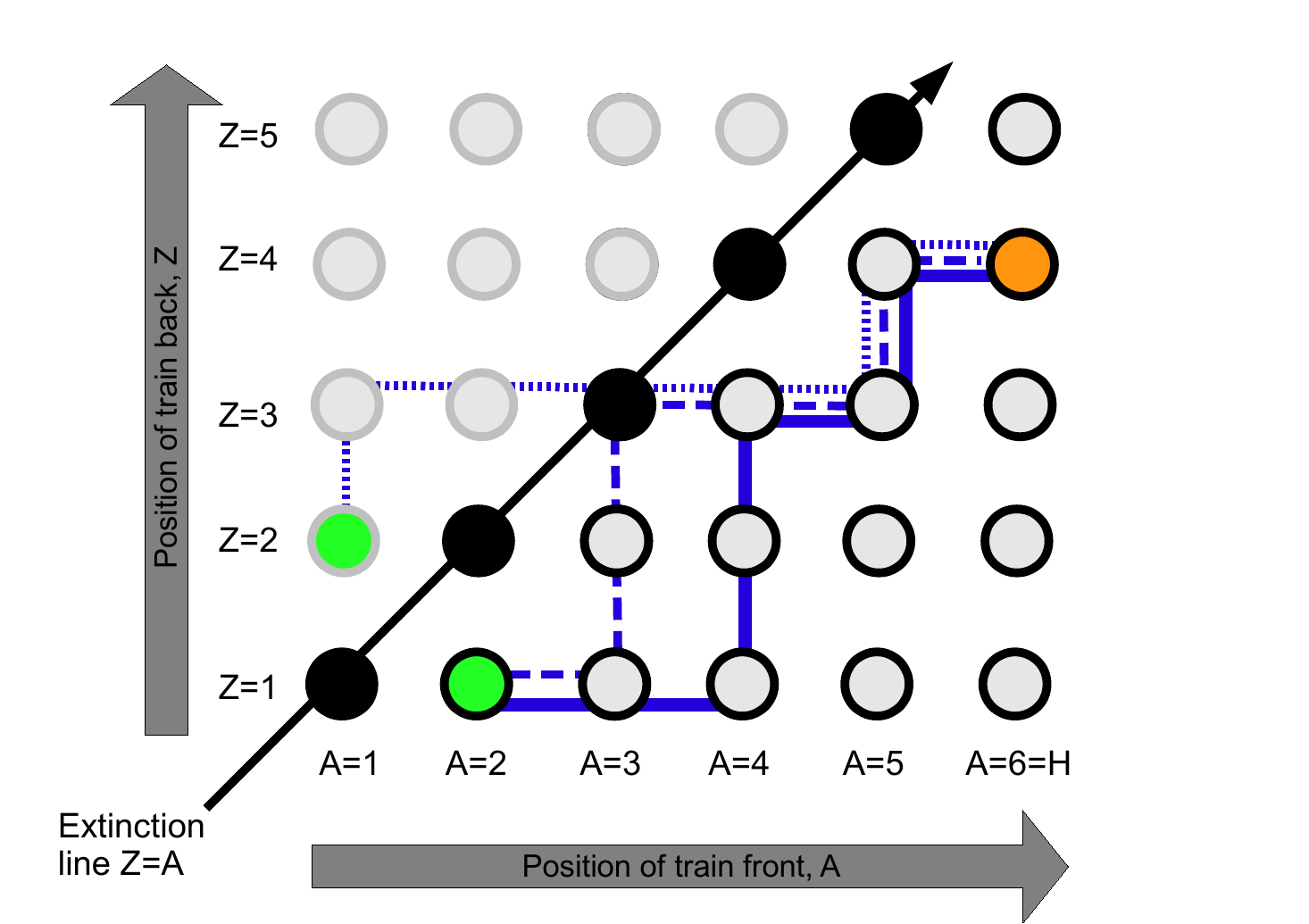}
\caption{Grid showing collection of possible train states. Permitted states (black outline), ghost states (grey outline) and extinction states (black fill) as well as a number of possible paths from our initial state (black outline, green fill) to a sample end state (black outline orange fill). Depicted are a permitted path (continuous), an invalid path (leading to extinction, long dash) and the associated ``ghost path'' from the reflection of our initial state to the sample end state (fine dash).}
\label{fig:reflection}
\end{figure}

In order to determine the expected train length, $T$, we consider the above process on a grid, where the horizontal axis represents the position of the front, $A_t$, and the vertical axis the back of the train, $Z_t$. Each point on the grid and below the diagonal, $A_t=Z_t$, represents a possible configuration of a train in the stem, see \fig{reflection}. 
All other points represent invalid configurations, which we refer to as ghost states. For each train, the initial configuration is $(A_0,Z_0)=(2,1)$, that is, the second stem node is a mutant, while the state of the first stem node is a resident (due to replacement on the fast time scale, leading to a slight underestimate).

Each train produces a trajectory or path on the grid that originates in $(A_0,Z_0)$ and ends at time $\tau$ once the train has reached its destination: the bottom of the stem, where $A_\tau=H$. However, if at any point in time $A_t\leq Z_t$ then this represents an \emph{invalid} path because the train has vanished. Every invalid path touches or crosses the diagonal $A_t=Z_t$ at least once. For a \emph{valid} path $A_t>Z_t$ must hold at all times. The expected train length, $T$, is the weighted average over all paths, with invalid paths being considered as having length zero. The number of valid paths can be calculated using the reflection principle \citep[p.~88]{koralov:book:2007}, which states that for every \emph{invalid} path a ghost path exists, starting from $(Z_0,A_0)$. The trajectory of a ghost path is the reflection of the corresponding invalid path along the diagonal $A_t=Z_t$, up to the point where the invalid path touches or crosses $A_t=Z_t$ for the first time. From then on the ghost path and the remainder of the invalid path coincide, see \fig{reflection}. 
In order to calculate the expected train length, $T$, we consider the train lengths based on \emph{all} paths and subtract all \emph{invalid} paths, to obtain the train length based on \emph{valid} paths only. The number of ghost paths corresponds to all paths starting from $(Z_0,A_0)$, the reflection of $(A_0,Z_0)$ and hence the name of the method. Having counted the number of paths we then weigh the corresponding train length by the probability of each possible path and obtain:
\begin{align}
\label{eq:TrainSum}
T &= (1-\alpha)^{H-2} \sum_{z=1}^{H-1} (H-z)\alpha^{z-1} \left[ \binom{H-3+z-1}{z-1}-\binom{H-3+z-1}{z-2} \right]
\end{align}
with $\alpha= 1/(1+r)$. We assume beneficial mutations, $r>1$, such that $0<\alpha<1/2$. All paths require $H-2$ steps that increase $A_t$ from the starting point at $2$ to the end point at $H$, which occurs with probability $1-\alpha$ for each step. The combinatorial sum then accounts for all possibilities and probabilities that $Z_t$ is increased along valid paths. In particular, the index variable $z$ indicates the position of the tail of the train and hence $H-z$ specifies the train length. The tail starts at $1$ and, for any valid path, reaches at most $H-1$. Because we are interested in the length of the train at the moment of arrival, the final step must be an increment of $A_t$. In particular, it follows that $z=H-1$ has zero valid paths -- a reassuring result as we know that no train could possibly have length one at the moment of its arrival. 
Note that we have used the convention that $\binom{n}{k}=0$ for $k<0$, which applies only if the tail remains at $Z_t=1$ and admits only a single valid path.

\subsection{Simplifying $T$}
\label{app:simpleT}
We now resort to  algebraic manipulation. Various binomial coefficient identities are used throughout:
\begin{align*}
T & =\ (1- \alpha)^{H-2} \sum_ {z=1}^{H-1} (H-z) \alpha^{z-1}  \left[ \binom {H+z-4}{z-1}-\binom {H+z-4}{z-2}   \right]\\
\intertext{\scriptsize using Pascal's rule}
 & =\ (1- \alpha)^{H-2} \sum_ {z=1} ^{H-1} (H-z) \alpha^{z-1} \left[  \binom {H+z-3}{z-1}-2\binom {H+z-4}{z-2} \right]\\
\intertext{\scriptsize splitting sum}
 & =\ (1- \alpha)^{H-2} \sum_ {z=1}^{H-1} (H-z) \alpha^{z-1}  \binom{H+z-3}{z-1}   -2 \alpha (1- \alpha)^{H-2} \sum_{z=1}^{H-1} (H-z) \alpha^{z-2} \binom {H+z-4}{z-2}\\
\intertext{\scriptsize changing second summation to obtain lower bound}
 & \ge\ (1- \alpha)^{H-2} \sum _ {z=1} ^{H-1} (H-z) \alpha^{z-1}  \binom{H+z-3}{z-1}  -2 \alpha (1- \alpha)^{H-2} \sum_{z=2}^{H} (H+1-z) \alpha^{z-2} \binom {H+z-4}{z-2}\\
\intertext{\scriptsize merging sums and relabelling indices}
 & =\ (1-2\alpha) (1- \alpha)^{H-2} \sum _ {z=0} ^{H-2} (H-1-z) \alpha^{z}  \binom{H+z-2}{z}\\
\intertext{\scriptsize expanding factor}
 & =\ (1-2\alpha) (1- \alpha)^{H-2} \sum _ {z=0} ^{H-2} (2(H-1) -(H-1)-z) \alpha^{z}  \binom{H+z-2}{z}\\
\intertext{\scriptsize using the combinatorial identity $(n+k)\binom{n+k-1}{k} = (n+k)\binom{n+k-1}{n-1}=n \binom{n+k}{n}= n \binom{n+k}{k}$}
 & =\ (H-1) (1-2\alpha) (1- \alpha)^{H-2} \sum _ {z=0} ^{H-2}  \alpha^{z}  \left[ 2 \binom{H+z-2}{z}- \binom{H+z-1}{z}  \right]\\
\intertext{\scriptsize extending the sum to $\infty$ can only decrease the lower bound because $2 \binom{n+k}{k}-\binom{n+k+1}{k}= \binom{n+k}{k} - \binom{n+k}{k-1} \le 0$ for $k>n$ and $(1-2 \alpha)>0$}
 & =\ (H-1) (1-2\alpha) (1- \alpha)^{H-2} \sum _ {z=0} ^{\infty}  \alpha^{z}  \left[ 2 \binom{H+z-2}{z}- \binom{H+z-1}{z}  \right] \\
\intertext{\scriptsize using $(1-\alpha)^{-n-1} = \sum _ {k=0} ^{\infty}  \alpha^{k} \binom{n+k}{k}$ since $|\alpha|<1$}
&= \ (H-1) (1-2\alpha) (1- \alpha)^{H-2} \left[2(1-\alpha)^{1-H} -(1-\alpha)^{-H}   \right]= (H-1) \left(\frac{1-2\alpha}{1-\alpha} \right)^2.
\end{align*}  
And so we have
\begin{align}
 T \ge (H-1) \left(1-\frac1r\right)^2.
\label{eq:TboundLow}
\end{align}
Hence, for $r>1$, the expected train length, $T$, can be made arbitrarily long by choosing a suitably long chain, $H$.

\subsection{Train collisions}
\label{app:collidingT}
The above derivation of the expected train length neglects the possibility that two trains may collide and merge, which introduces a source of error. When two trains collide, the first finds itself being erased at greater than the expected rate, effectively reducing its train length. Therefore, collisions \emph{decrease} the expected train length $T$ -- despite the fact that merging trains may lead to longer overall lengths -- hence \eq{TrainSum} overestimates the expected train length. An upper bound on our true $T$ is found by assuming that trains never interfere with one another, and simply using \eq{TrainSum}. Conversely, a lower bound for $T$ is obtained by assuming that the second train completely eradicates the first train. In order to find this lower bound we need to determine the probability that a train collision occurs. An upper bound on this collision probability is given by the probability that a second train is generated while another train is still occupying the stem. This can be formulated in terms of a negative binomial distribution where the generation of a new train counts as a ``success'' while a decrease in length of the existing train in the stem counts as a ``failure''. In each time step a new train is generated with probability $r^2/(F_t (L+r-1))$ whereas the probability that the existing train length decreases, i.e. the resident directly behind the train reproduces, is at least $1/F_t$ (exactly $1/F_t$ along the stem, but greater for the first stem node). After $H$ failure events we know that the stem must be cleared and contain only residents. Therefore, train collisions occur at most with the probability that a new train is generated prior to $H$ failure events:
\begin{align}
\label{eq:no2trains}
P(\text{no $2^{nd}$ train})= \left(1-\frac{r^2}{L+r-1+r^2}\right)^H > 1-H \frac{r^2}{L+r-1+r^2}.
\end{align}
The inequality in \eq{no2trains} results from expanding and then truncating the power term. This is permitted, because when expanding the alternating sum for sufficiently large $L$, the absolute value of subsequent terms is strictly decreasing. Thus, the chance that a second train is launched while another one still occupies the stem is at most 
\begin{align}
\label{eq:ep3}
\epsilon_3 = \frac{H r^2 }{L+r-1+r^2}
\end{align} 
 and becomes small for $L\gg H r^2$. Thus, the true expected train length lies between $T$ and $T(1-\epsilon_3)$.

\section{Interaction of trains with root node}
\label{app:stemBase}
Here we calculate the probability that a train of mutants, which arrives at the base of the stem with an initial length $l$, succeeds in taking over the root node and placing a new mutant in one of the reservoirs. We adapt the technique used in \citet{diaz:PRSA:2013}, considering a finite state Markov process with two absorbing states: either a new mutant is placed in one reservoir, or the mutant train has disappeared. All other states represent a particular train length with the root node occupied by either resident or mutant.

For each state, the probability of eventually succeeding is $p_{i \updownarrow}$, where $i\leq l$ indicates the current train length and $\scriptstyle\updownarrow$ indicates whether the root node is occupied by a mutant, $\scriptstyle\uparrow$, or resident, $\scriptstyle\downarrow$. Clearly $p_{0 \downarrow}=0$ because the train has disappeared, the root is a resident and hence an absorbing state has been reached. Similarly, $p_{0 \uparrow}=r/(B+r)$, denotes the probability that the mutant in the root node reproduces before being replaced by the offspring of residents in any of the $B$ branches.
By examining all possible transitions we obtain: 
\begin{align}
(B+r ) p_{i \uparrow}  =&\ (B-1)p_{i \downarrow} + r + p_{i-1 \uparrow}  \\
(r+1 ) p_{i \downarrow}  =&\  r p_{i \uparrow}+  p_{i-1 \downarrow}
\end{align}
If the train is $i$ mutants long, and the root is a mutant, a ``success'' occurs with relative weight $r$, whereas the root mutant is lost, with relative weight $B-1$ (because there are $B-1$ other branches), or our train may erode, with relative weight $1$, leaving the root node unchanged. If the root node is a resident, then the only possible actions are the replacement of the root node, or eroding the train  from behind, with relatively probabilities $r$ and $1$ respectively. Finally, we have coefficients on the left hand side to normalize over all possible courses of action. Written as a matrix equation this gives:
\begin{align*}
\left[ \begin{array}{cc} \
B+r & 1-B \\ 
-r & r+1 \end{array}\right]
\left[ \begin{array}{c}
p_{i \uparrow}\\ p_{i \downarrow}
\end{array} \right]
 =&\ 
\left( \left[ \begin{array}{c}
r\\ 0
\end{array} \right]
 +
\left[ \begin{array}{c}
p_{i-1 \uparrow}\\ p_{i-1 \downarrow}
\end{array} \right]
 \right)
\end{align*}
which yields
\begin{align}
\left[ \begin{array}{c}
p_{i \uparrow}\\ p_{i \downarrow}
\end{array} \right]
 =&\ \frac{1}{B +2r +r^2}  \left[ \begin{array}{cc} \
r+1 & B-1 \\ 
r & r+B \end{array}\right]
\left( \left[ \begin{array}{c}
r\\ 0
\end{array} \right]
 +
\left[ \begin{array}{c}
p_{i-1 \uparrow}\\ p_{i-1 \downarrow}
\end{array} \right]
 \right).
\label{eq:PUpDownEquations}
\end{align}
We now calculate both upper and lower bounds on the expected success probability upon arrival, $\mathbb{E} (p_{l \downarrow})$. Let us start with the upper bound. 

By neglecting terms in the normalization term and increasing several of the matrix entries, we find an upper bound on the R.H.S. of \eq{PUpDownEquations}. This works because $p_{i \uparrow}, p_{i \downarrow} \ge 0$, and thus we are making the R.H.S. strictly more positive. It also significantly simplifies our equation to
\begin{align*}
\left[ \begin{array}{c}
p_{i \uparrow}\\ p_{i \downarrow}
\end{array} \right]
 < \frac{1}{B+2r +1}   \left[ \begin{array}{cc} \
r+1 & B+r \\ 
r+1 & r+B \end{array}\right]
\left( \left[ \begin{array}{c}
r\\ 0
\end{array} \right]
 +
\left[ \begin{array}{c}
p_{i-1 \uparrow} \\ p_{i-1 \downarrow}
\end{array} \right]
 \right)
\end{align*}
where the inequality applies element-wise. Substituting $p_{0 \uparrow} = r/(B+r) , p_{0 \downarrow}=0$ into the above gives
\begin{align*}
p_{1 \updownarrow} <\frac{r^2 +r}{B+2r+1} \left(1+ \frac{1}{B+r} \right)
\end{align*}
Using the fact that the upper bounds for $p_{i \uparrow}$ and $p_{i \downarrow}$ are equal, along with the fact that the denominator of the fraction is equal to the row sum of the transition matrix gives:
\begin{align*}
p_{i \updownarrow} <\frac{r^2 +r}{B+2r+1} + \frac{B+2r+1}{B+2r+1}p_{i-1 \updownarrow}.  \\
\intertext{By induction we find}
 p_{i \updownarrow} <\frac{r^2 +r}{B+2r+1} \left(i  + \frac{1}{B+r} \right).
\end{align*}
This yields an upper bound for $p_{i \uparrow}$, which we then use to calculate a tighter bound for $p_{i \downarrow}$.
From the last line of \eq{PUpDownEquations} we derive:
\begin{align*}
p_{i \downarrow} < \frac{r}{B} (r+p_{i-1 \uparrow}) +p_{i-1 \downarrow}
\end{align*}
which leads to
\begin{align*}
p_{i \downarrow} < \sum\limits_{n=1}^{n=i} \frac{r}{B} (r+p_{n-1 \uparrow}) \le \frac{i r^2}{B} \left( 1 + \frac{r +1}{(B+2r +1)(B+r) }+  \frac{H-1}{2} \frac{r +1}{B+2r+1} \right).
\end{align*}
Thus we end up with $p_{i \downarrow} < i r^2(1 + \epsilon_{4+})/B$, where
\begin{align}
\label{eq:ep4Plus}
 \epsilon_{4+}= \frac{1+r}{(B+2r+1)(B+r)} + \frac{(H-1)(r+1)}{2B+4r+2} 
\end{align}
 is our error term. This error term can be made arbitrarily small for sufficiently large $B$.

For the lower bound, we must deal with the possibility that a small number of other branches, $\delta$, also contain mutants.  
Because the above train collision argument (\app{collidingT}) is based on at most a single mutant existing in the reservoir of each branch, we wish to only consider reproductive events, which place new mutants into branches that currently contain no mutants and neglect the rest (which we can do, as this calculation desires a \emph{lower} bound). Further, we must contend with the fact that trains from other mutants may compete for control of the root node. Although generically we do not expect to encounter other trains, we calculate our lower bound as if all other mutant occupied branches have mutants at the base of their stems at all times. This arrangement, while unrealistic, describes the situation which minimizes the success probability of a given train, and is thus useful for finding lower bounds. The following equations are written under the assumption of this worst case scenario (worst from the perspective of the train we are focusing on):
\begin{align*}
\left[ \begin{array}{cc} \
B+r+\delta(r-1) & 1-B-\delta(r-1) \\ 
-r & r+1 \end{array}\right]
\left[ \begin{array}{c}
p_{i \uparrow}\\ p_{i \downarrow}
\end{array} \right]
 =
\left( \left[ \begin{array}{c}
r \frac{B-\delta}{B} \\ 0
\end{array} \right]
 +
\left[ \begin{array}{c}
p_{i-1 \uparrow}\\ p_{i-1 \downarrow}
\end{array} \right]
 \right)
\end{align*}
which can be rewritten as
\begin{align*}
\left[ \begin{array}{c}
p_{i \uparrow}\\ p_{i \downarrow}
\end{array} \right]
 = \frac{1}{B+\gamma(r-1)+2r +r^2}   \left[ \begin{array}{cc} \
r+1 & B+\delta(r-1)-1 \\ 
r & r+B+\delta(r-1) \end{array}\right]
\left( \left[ \begin{array}{c}
r\frac{B-\delta}{B}\\ 0
\end{array} \right]
 +
\left[ \begin{array}{c}
p_{i-1 \uparrow}\\ p_{i-1 \downarrow}
\end{array} \right]
 \right).
\end{align*}
By  ignoring the positive effects of $ p_{i-1 \uparrow}$ on $ p_{i \downarrow}$ we form the inequality
\begin{align*}
p_{i \downarrow}
 > \frac{1}{B+\delta(r-1)+2r +r^2}  
\left(  r^2 \frac{B-\delta}{B}
 +
(r+B+\delta(r-1)) p_{i-1 \downarrow}
 \right)
\end{align*}
for $i \ge 1$. By induction this leads to
\begin{align*}
p_{i \downarrow}
 &> \frac{B-\delta}{B} \frac{r^2}{B + \delta(r-1)+2r +r^2}  
  \sum\limits_{n=0}^{i-1} \left( \frac{B+r+ \delta(r-1)}{B+ \delta(r-1)+r^2+2r} \right)^n\\
 &= \frac{B-\delta}{B} \frac{r^2}{B + \delta(r-1)+2r +r^2}  
 \frac{ 1-\left( \dfrac{B + \delta(r-1)+r}{B+ \delta(r-1)+r^2+2r} \right)^i }{ 1-\dfrac{B + \delta(r-1)+r}{B+ \delta(r-1)+r^2+2r} }
\end{align*}
and finally to
\begin{align*}
p_{i \downarrow}
 >  \frac{B-\delta}{B} \frac{r^2}{r +r^2}  
\left( 1-\left(1- \frac{r^2+r}{B + \delta(r-1)+r^2+2r} \right)^i \right).
\end{align*}
As long as $i (r^2+r)/(B + \delta(r-1)+r^2+2r) < 1$, the series expansion of the inner bracket gives an alternating sequence with monotone decreasing absolute terms. Therefore, we can truncate the series after three terms while still preserving the inequality because the sum of the first three terms is greater than any subsequent sum.
This leads to: 
\begin{align*}
p_{i \downarrow}>\frac{B-\delta}{B} \frac{ir^2}{B + \delta(r-1)+r^2+2r}
\left( 1 -   \frac{i-1}{2}  \frac{r^2+r}{B + \delta(r-1)+r^2+2r} \right)
\end{align*}
If we rearrange the above, remembering that $1/(1+x)<1-x$ whenever $x>-1$, we find
\begin{align*}
p_{i \downarrow}> \frac{i r^2}{B} \left(1- \frac{\delta}{B} - \frac{\delta(r-1)+r^2 +2r}{B}- \frac{H-1}{2} \frac{r^2+r}{B} + O(B^{-2}) \right).
\end{align*}
We are free to drop the very small positive terms at the end, and find the result
$p_{i \downarrow} > (1- \epsilon_{4-}) i r^2/B $, where the error term
\begin{align}
\label{eq:ep4Minus}
\epsilon_{4-}= \frac{2\delta r+r^2 +3r+ H(r^2+r)}{2B}
\end{align}
can be made small whenever $\delta,H \ll B$. Thus $\mathbb{E}(p_{l \downarrow}) >  T r^2 B^{-1}(1 - \epsilon_{4-})$. Armed with the constrains $  (1- \epsilon_{4-}) l r^2/B < p_{l \downarrow} < (1 + \epsilon_{4+}) l r^2/B$, we note that  $\mathbb{E}(p_{l \downarrow}) \approx \mathbb{E}(l r^2 /B) = T r^2 / B$.

\section{Bounds on fixation probabilities}
\label{app:fixBound}

\subsection{Upper bound}
\label{app:upperFixBound}
Starting with a single mutant in one reservoir initially, we first must have two before we can have $BL$ reservoir mutants. Thus, the probability of transitioning from $1$ to $2$ mutants in the reservoir serves as a straight forward upper bound for the mutants fixation probability.
To further simplify the upper bound, we make several optimistic (from the mutants point of view) assumptions: (i) the original mutant appears in a reservoir node (ignoring $\epsilon_0$, see \eq{ep0}); (ii) We slightly increase the train launch probability, dividing by $L$ rather than $L+r-1$ (ignoring $\epsilon_1$, see \eq{ep1}) (iii) no detrimental effects based on our initial conditions (ignoring $\epsilon_2$, see \eq{ep2}); (iv) no train collisions (ignoring $\epsilon_3$, see \app{collidingT}); and finally (v) we use the upper bound for the probability that a train succeeds in producing another reservoir mutant ($(1+\epsilon_{4+}) T r^2/B$, see \eq{ep4Plus}). 

A single mutant in any reservoir produces a new train with a probability of at most $r^2/(F_t L)$ per time step. Subsequently, each train succeeds in placing another mutant in any reservoir with a probability of at most $(1+\epsilon_{4+})T r^2/B$. At the same time, the root node has a probability of at least 
\begin{align*}
\frac{B-1}{B+r-1} \frac{1}{F_t B L}
\end{align*}
to remove the mutant node from the reservoir.
Thus, the chance of the mutant producing a successful train before being erased by the root node is \emph{at most} $\rho_{H+}$ with
\begin{align}
\label{eq:UpperBound}
\rho_{H} \le \rho_{H+} =\frac{T r^4(1+\epsilon_{4+}) }{T r^4(1+\epsilon_{4+}) + \frac{B-1}{B+r-1}} \approx 1- \frac{1}{T r^4 +1}.
\end{align}
The approximation in \eq{UpperBound} becomes exact in the limit of large $B$. Moreover, because $T<H$, we can replace $T$ by $H$ in \eq{UpperBound} to obtain a simpler and more generous upper bound.

\subsection{Lower bound}
\label{app:LowerFix}%
On the slow timescale the dynamics of the reservoir can be approximated by a random walk $X_t$ on the number of mutants in the reservoir. 
Consider the random walk on the integers from $0$ to $BL$ with forward bias $\gamma$ for $X_t< \delta\ll  B$ and no bias for $X_t \ge \delta$. Because the chance of any particular reservoir mutant replacing any particular reservoir resident is always higher than the converse, we can be sure that some forward bias persists for all $X_t$. Unfortunately, because the analytic arguments in the previous sections fail for many reservoir mutants, we must make the conservative assumption that no bias applies in this region. The fixation probability of this random walk acts as a lower bound on the fixation probability of the true process. In the following, we assume that $H$ and $r$ are fixed and that $B,L \gg H, \delta$. 

In order to determine the fixation probability of the random walk $X_t$, we construct a martingale, $Q(X_t)$.
A martingale is a function of a random variable such that the expected value of the martingale in the next time step is equal to the current value:
\begin{align}
\mathbb{E}(Q(X_{t+1}) |Q(X_t) )= Q(X_t).
\end{align}
For $Q(X_t)$ to be a martingale, we require:
\begin{align}
Q(k) =&\ 
\begin{cases}
\begin{array}{lll}
\dfrac{\gamma}{1+ \gamma } Q(k+1) + \dfrac{1}{1+ \gamma } Q(k-1) & \ 0<k<\delta & \ \text{(forward bias)}\\[10pt]
\dfrac{1}{2} Q(k+1) + \dfrac{1}{2} Q(k-1) & \ \delta \le k < BL & \ \text{(no bias)} .
\end{array}
\end{cases}
\end{align}
These constraints admit the solution $Q(k)= \gamma^{-k}$ for $k<\delta$, and $Q(k)=A k +D$ for $k \ge \delta$. For $Q(k)$ to satisfy the martingale conditions as needed, we demand $\delta \in \mathbb{N}$.
The constants $A,D$ are determined by connecting the solutions for the two regions. In particular, 
\begin{align*}
\gamma^{-\delta} =&\ Q(\delta)= A \delta +D
\end{align*}
must hold such that $Q(\delta)$ is well defined and 
\begin{align*}
2 (A \delta +D)=&\ \gamma^{-\delta+1} +  A \delta+A +D
\end{align*}
to satisfy the martingale property at $\delta$. Thus, 
\begin{align*}
A=&\ \gamma^{-\delta}(1-\gamma) \\
D=&\ \gamma^{-\delta}(1-\delta(1-\gamma))
\end{align*}
 which yields 
\begin{align*}
Q(0)=&\ 1 \\
Q(1)=&\ \gamma^{-1} \\ 
Q(BL)=&\ \gamma^{-\delta} + \gamma^{-\delta}(1-\gamma) (BL-\delta).
\end{align*}

Let $\tau$ be the first time we reach one end of our random walk. Because $Q(k)$ is bounded for all relevant values of $k$ we are able to invoke the optional stopping theorem \citep[p. ~210]{klenke:book:2006}, which renders
\begin{align*}
Q(1)=Q(X_0)= \mathbb{E}(Q(X_\tau))= Q(0) P(0) + Q(BL) P(BL),
\end{align*}
where $P(0)$ and $P(BL)$ represent the probabilities of reaching either end of our random walk. 
Using $P(0) =1- P(BL)$ we find
\begin{align}
P(BL) = \frac{Q(1)-Q(0)}{Q(BL) -Q(0)}= \frac{1- \gamma^{-1} }{1-\gamma^{-\delta} -\gamma^{-\delta}(1-\gamma)(BL-\delta) }.
\end{align}
In order to keep error terms small, we must select $\delta$ such that $\gamma^\delta \gg BL$ and $\delta \ll B,L$. 
As long as $B$ and $L$ are large and sufficiently similar, this is possible provided $\gamma>1$ (see below). Thus, for any choice of $H$, and for any $r>1$ we can select $B$ and $L$ such that
\begin{align}
\label{eq:MartingFinal}
P(BL) =&\ \frac{Q(1)-Q(0)}{Q(BL) -Q(0)}= \frac{1- \gamma^{-1} }{1+ \epsilon_5 }\\
\intertext{with}
\epsilon_5 =&\ \gamma^{-\delta} \left( (\gamma-1)(BL-\delta) -1 \right) \ll 1.\notag
\end{align}
In order to find an upper bound on $\epsilon_5$, we require a lower bound on the forward bias $\gamma$. In particular, we would like to show that $\gamma>1$ and hence that $\epsilon_5$ can be made small. This can be seen by taking the lower bound on the production rate of successful trains, $(1-\epsilon_1)(1-\epsilon_3)(1-\epsilon_{4-}) Tr^4/(B L F_t)$), and comparing to our upper bound on the removal probability for reservoir mutants, $1/(B L F_t)$. 
This represents the eventual forward bias after the top of the stem has been replaced at least once. To account for the possibility of mutant loss before the top of the stem has been replaced we must consider $\epsilon_2$, which acts as an additive penalty (because it only applies once per reservoir mutant). This yields
\begin{align}
\label{eq:gamma}
\gamma\geq r^4 T (1-\epsilon_1)(1-\epsilon_3)(1-\epsilon_{4-}) -\epsilon_2.
\end{align}

In the limit of large $B$ and $L$ all error terms tend to zero. Thus, to show $\gamma>1$, it is sufficient to show that $r^4 T>1$. 
Recalling that $A_t-Z_t$ represents the length of a train at time $t$ (see \app{Mirror}), and noting that it is submartingale 
(the expected future value is greater than the current value) whenever $r>1$, we can easily show that $T= \mathbb{E}(A_{\tau} -Z_{\tau}) \ge A_{0} -Z_{0}=1$, and thus, in the limit, $\gamma \ge r^4>1$. Thus $\epsilon_5$ can indeed be made arbitrarily small.

Substituting \eq{gamma}, the lower bound of $\gamma$ into \eq{MartingFinal} yields a lower bound on $P(BL)$, which in turn provides a lower bound on the fixation probability, see \eq{rhoH}. To simplify the lower bound substitute the lower bound for $T$, $(H-1)(1-1/r)^2)$, into \eq{gamma} and obtain a looser bound on fixation probability.

\section{Bringing it all together}%
\label{app:summary}%
Collating \eq{MartingFinal} and \eq{UpperBound}, we find:
\begin{align}
\label{eq:FINALBOUNDS}
& \frac{1-\epsilon_0 }{1+ \epsilon_5 }\left(1- \frac{1}{r^4 T(1-\epsilon_1)(1-\epsilon_3)(1-\epsilon_{4-}) -\epsilon_2}   \right) \le \rho_H \le 1- \frac{B-1}{(B+r-1)T r^4(1+\epsilon_{4+}) + B-1}
\end{align}
with the train length, \eq{TrainSum}
\begin{align*}
T =&\ (1-\alpha)^{H-2} \sum_{z=1}^{H-1} (H-z)\alpha^{z-1} \left[ \binom{H-3+z-1}{z-1}-\binom{H-3+z-1}{z-2} \right],  \\
 &\ (H-1)(1-r^{-1})^2 \le T \le H,  \\ 
\intertext{the chance that the initial mutant is not in the reservoir, \eq{ep0}}
\epsilon_0 =&\ \frac{1+HB}{BL+1+HB}, \\ 
\intertext{Simplifying approximation in train launch probability, \eq{ep1}}
\epsilon_1= \frac{r-1}{L+r-1} \\
\intertext{the chance that the initial mutant is removed before it reproduces, \eq{ep2}}
\epsilon_2 =&\ \frac{1}{1+ BL^2}, \\ 
\intertext{the chance of train collisions, \eq{ep3}}
\epsilon_3 =&\  H \frac{r^2}{L+r-1+r^2}, \\ 
\intertext{the lower bound for train success, \eq{ep4Minus}}
\epsilon_{4-}=&\ \frac{2\delta r+r^2 +3r+ H(r^2+r)}{2B}, \\ 
\intertext{the upper bound for train success, \eq{ep4Plus}}
\epsilon_{4+} =&\  \frac{ (r +1)}{B+2r +1} \left( \frac{1}{B+1} + \frac{H-1}{2} \right), \\ 
\intertext{and finally the Martingale error term, \eq{ep4Plus}}
\epsilon_5 =&\ \gamma^{-\delta} [ (\gamma-1)(BL-\delta) - 1]. 
\end{align*}
Because many of the error terms are dependent on both $B$ and $L$ it makes sense to take the limit of both simultaneously. This will force all error terms to zero as long as  $B,L \gg H,\delta$ and $\epsilon_5 \to 0$.   In particular, in the limit $B \to \infty$, with $\sqrt{B}-1 < \delta \le \sqrt{B}$ and $L=B$ all error terms tend to zero. Other variations, such as $L=B^\beta$, $\beta >0$ are possible, although it is suspected that relations of the form $L= \beta^B$ would prove problematic, as we would then need to reconcile the bounds $\delta \ll B$ and $\gamma^\delta \gg B \beta^B$, a problem we do not run into for $L=B^\beta$.

Simpler and looser upper and lower bounds independent of $T$ are obtained by substituting upper and lower bounds for $T$ (respectively) into \eq{FINALBOUNDS}.

\section{Deleterious mutations, $r<1$}
\label{app:briefNote}
As a final remark, we note that the above arguments are based on the assumption of a beneficial mutation, $r>1$. In addition to promoting beneficial mutations, an evolutionary amplifier must also suppress the fixation of deleterious mutations, $r<1$. Here we argue that a deleterious mutant indeed disappears almost surely on sufficiently large superstars.

Consider a single mutant with fitness $1$, in a population of residents with fitness $1/r$. Note that we can rescale fitness without changing the dynamics of the system, since fitness is never used in an absolute sense, only relative fitness matters. We next observe that all calculations performed previously with respect to a rare mutant with a fitness advantage would now apply to a resident, if it were to become rare -- that is, if residents were rare, we would expect trains of residents to propagate down the stem (incrementing with probability $1/(r+1)> 1/2$ and shrinking with probability $\alpha= r/(1+r)<1/2$). The same martingale argument that we previously used to find a lower bound on mutant fixation probability can now be used to obtain a lower bound on resident fixation probability. 
This time we use the random walk $X_t$ to track the number of residents in our reservoir nodes. Thus $X_0=BL-1$.
This leads to a formula for the fixation probability of \emph{residents}
\begin{align}
1-\rho_H \ge&\ P(BL)= \frac{Q(0)-Q(BL-1)}{Q(0) -Q(BL)} \notag\\
=&\ \frac{1-\gamma^{-\delta} -\gamma^{-\delta}(1-\gamma)(BL-1-\delta) }{1-\gamma^{-\delta} -\gamma^{-\delta}(1-\gamma)(BL-\delta) } \notag\\
=&\ 1 -\frac{  \gamma^{-\delta}(\gamma-1) }{1-\gamma^{-\delta} -\gamma^{-\delta}(1-\gamma)(BL-\delta) } \notag\\
\intertext{and hence}
\rho_H \le&\ \frac{  \gamma^{-\delta}(\gamma-1) }{1-\gamma^{-\delta} -\gamma^{-\delta}(1-\gamma)(BL-\delta) } \notag\\
 \le&\ \gamma^{1-\delta}.
\label{eq:detrmut}
\end{align}
In the above we require $\gamma$ to be the bias in favour of the \emph{resident}. In order to calculate it we need to find the expected train lengths of resident trains. Structurally the train equation is the same as previously, but with $r$ replaced by $1/r$ (the fitness of  residents). Because $1/r$ is now greater than one, all train length arguments that previously depended on $r>1$ for beneficial mutations can now be applied to the resident with $1/r>1$. Thus, we can show that the expected length of resident trains is large. Because $\gamma \approx  r^{-4}T$, large $T$ leads directly to a large $\gamma$ (in the residents favour), and $\gamma^{-\delta} BL \ll 1$. The Martingale argument will significantly underestimate the fixation probability for the resident and thus overestimate the probability of mutant fixation. 
In order to derive exact bounds on the new bias in favour of residents, $\gamma$, we must apply error terms similar to $\epsilon_2, \epsilon_3 , \epsilon_{4-}$, to obtain lower bounds on the effectiveness of ``resident trains''. No error term equivalent to $\epsilon_0$ arises because the possibility that the initial mutant is not placed in a reservoir node only reduces the mutant fixation probability. 
Thus, the fixation probability of a deleterious mutant is bounded by \eq{detrmut} and can be made arbitrarily small for sufficiently big $H, L$ and $B$.

\end{document}